
\def\cl{\centerline}
\def\md{\medskip}
 \def\ds{\displaystyle}
 \def\sqr#1#2{{\vcenter{\vbox{\hrule height.#2pt
               \hbox {\vrule width.#2pt height#1pt \kern#1pt
                     \vrule width.#2pt}
               \hrule height.#2pt}}}}
\def\square{{\mathchoice\sqr57\sqr57\sqr{2.1}3\sqr{1.5}3}}
\def\proclaim #1. #2\par{\medbreak
     {\smps #1.\enspace}{\it #2}\par\medbreak}

\def\cB{{\cal B}}
\def\cC{{\cal C}}
\def\cD{{\cal D}}
\def\cG{{\cal G}}
\def\cL{{\cal L}}

\def\cS{{\cal SP}}
\def\cX{{\cal X}}

\def\tKG{{\widetilde { K\cG}}}

\def\kplus{{\widetilde { K\cG}_+(\chi)}}
\def\bns{[1]}
\def\br{[2]}
\def\bs{[3]}

\def\brown{[4]}
\def\cf{[5]}
\def\droms{[6]}
\def\d2{[7]}
\def\cehpt{[8]}
\def\gromov{[9]}
\def\hm{[10]}
\def\topgraph{[11]}
\def\sds{[12]}
\def\stallings{[13]}
\def\vw{[14]}
\def\w{[15]}
 3
 2
\font\small = cmr9
\font\sbaps = cmcsc10 scaled\magstep 2
\font\scaps = cmcsc10 scaled\magstep 1
\font\smps = cmcsc10
\tolerance=1000
\narrower
\magnification = 1200
\baselineskip = 13pt
{
\baselineskip = 14pt
\cl {\sbaps THE BIERI-NEUMANN-STREBEL INVARIANTS}
\md
\cl {\sbaps FOR GRAPH GROUPS}
\bigskip
}
\cl {{\scaps JOHN MEIER}~{\smps and}~{\scaps LEONARD VANWYK}}
\bigskip
\
\md
{
\small
\baselineskip = 10pt
\cl {\smps Abstract}
\bigskip
Given a finite simplicial graph $\cG$,
the graph group $G\cG$ is the group with generators in
one-to-one correspondence with the vertices of $\cG$ and with
relations stating two generators commute if their associated
vertices are adjacent in $\cG$.   The Bieri-Neumann-Strebel
invariant can be explicitly described in terms of the original
graph $\cG$ and hence there is an explicit
description of the distribution
of finitely generated normal subgroups of $G\cG$
with abelian quotient.  We construct Eilenberg-MacLane
spaces for graph groups and find partial extensions of this
work to the higher dimensional invariants.
\smallskip
}
\
\bigskip
\cl {\it Introduction}
\bigskip
Let $G$ denote a finitely generated group.\footnote {}
{1991 {\it Mathematics Subject Classification, Primary:
{\rm 20E07}, Secondary: {\rm 20J05}.}}
In several papers a collection of geometric invariants $\Sigma^k(G)$
($k$ a positive integer) was
developed, each of which is a subset of the character sphere
for $G$.  These  Bieri-Neumann-Strebel invariants have proven to be
quite rich.  For instance if $G \simeq \pi_1(M)$
for a smooth compact manifold $M$, then $\Sigma^1(G)$
yields information of the existence of circle fibrations of $M$,
and if in addition $M$ is a $3$-manifold then $\Sigma^1(G)$
can be described in terms of the Thurston norm.  If
$N$ is a normal subgroup of $G$ with $G/N$ abelian, then
whether $N$ has the finiteness property $F_k$ is
measured by
$\Sigma^k(G)$.  In particular $\Sigma^1(G)$ measures
the finite generation of normal subgroups
of $G$ with abelian
quotient.  (See \bns, \br, \bs\ and the
references cited there.  Further background will be
given in section 1.)  Regretably, \lq\lq $\Sigma^1(G)$ is,
in general, very difficult to compute." \bs
\md
Free partially commutative (FPC) monoids were first introduced by P.
Cartier and D. Foata \cf\ in order to study combinatorial problems
involving rearrangements of words.  In the last ten years,
they have been studied by both computer scientists and mathematicians.
The corresponding FPC groups are known as
graph groups.  A finite simplicial graph $\cG$ induces the following
presentation of a group $G\cG$:
$$
\langle V(\cG)\ \big |\ xy=yx \ \forall x,y\in V(\cG)\hbox{ such that }
x \hbox{ and } y \hbox{ are adjacent}
\rangle,
$$
where $V(\cG)$ denotes the vertex set of $\cG$.  A group $G$ is called
a {\it graph group\/} provided there exists some finite simplicial graph
$\cG$ such that $G\simeq G\cG$.  Given any graph group $G\cG$  we will often
assume the above presentation and we will identify the set of generators
of $G\cG$
with the vertex set $V(\cG)$.  Notice that if $\cG$ is completely
disconnected, then $G\cG$ is the free group on $V(\cG)$, while if
$\cG$ is a complete graph, then $G\cG$ is the free abelian group on
$V(\cG)$.  If $Y\subseteq V(\cG)$ and all $x,y\in Y$ are
adjacent, then $Y$
is called a {\it (commuting) clique\/}.
\par
Much is already known about graph groups.  For instance,
the word problem for graph groups was shown to be solvable by C.
Wrathall  \w.   Further, graph groups have been shown to satisfy
quadratic isoperimetric and linear isodiametric inequalities
\topgraph. The conjugacy problem was solved independently by
C. Wrathall and H.  Servatius.  C. Droms, H. Servatius, and B. Servatius
have studied various subgroups of graph groups (\droms, \d2, \sds).  Graph
groups were shown to be biautomatic by both authors: this
was proven by the first author
 in a joint paper with S.
Hermiller using the more general notion of a ``graph product'' of
groups \hm, and independently by the second author  \vw.
\md
Sections 2 through 5 of this paper
 establish
 a description of
the invariant $\Sigma^1(G\cG)$ in terms of the
defining graph $\cG$ (Theorem 5.1).  In particular we give an explicit
decomposition of $\Sigma^1(G\cG)$ into open simplices
in a standard simplicial decomposition of the $n$-sphere.
(The number $n$ depends on the number of vertices in $\cG$.)
\par
In section 6 we use this decomposition of $\Sigma^1(G\cG)$
to describe the finitely generated normal subgroups of
$G\cG$ with abelian quotient.  In particular,
let $\chi$ be a rational character of a graph group $G\cG$,
that is, $\chi$ is a map from $G\cG$ onto an infinite cyclic
group.  Let $\cL(\chi)$ be the full subgraph of $G\cG$
generated by vertices of $G\cG$ which do not map to $0$
under $\chi$.   Recall that a subgraph $\cL$ of
$\cG$ is {\it dominating} if for every vertex $v \in \cG\backslash\cL$,
$d(v,\cL)=1$.
\md
\proclaim Theorem 6.1.  Let $\chi$ be a rational character
of a graph group $G\cG$.  The kernel of $\chi$ is
finitely generated if and only if $\cL(\chi)$ is a
connected and dominating subgraph of $\cG$.
\hfill $\square$
\par
In section 7 we describe how to construct Eilenberg-MacLane
spaces for graph groups using  \lq\lq non-positive curvature"
techniques of Gromov.  In the final section we use these
Eilenberg-MacLane spaces to give partial descriptions of some
of the higher Bieri-Neumann-Strebel invariants.
\par
There is not currently a specific description of the
higher invariants for general graph groups.  However
we establish that in certain cases the invariants
have the same sort of stability as is present for
$1$-relator groups and three manifold groups.
\eject
\proclaim Theorem 8.3.  Let $G\cG$ be a graph group
which can be expressed
as the direct product of a non-trivial free abelian group and
a graph group with disconnected graph.
Then $\Sigma^1(G\cG) = \Sigma^k(G\cG)$ for all positive
integers $k$.\hfill $\square$
\par
It is known that Theorem 8.3 is
not true for arbitrary graph groups
\md
We thank Cliff Reiter for providing figures 1 and 2
in section 4 and Walter Neumann for an excellent
talk in the CUNY group theory seminar which
 inspired us to place this work within the framework
of the Bieri-Neumann-Strebel invariants.
Initial work on this project was done by the second
author as part of his doctoral dissertation under the direction
of the late Craig C. Squier.
\vskip 1cm
\cl {\it {\rm 1.} Bieri-Neumann-Strebel invariants}
\bigskip
The reader is advised to read \bns, \br\  or \bs; we will
only briefly review the main definitions and relevant
theorems.
\md
Any non-zero map from a
finitely generated group $G$ to the
additive group of the reals
is a {\it character} of $G$.
This is slightly non-standard
in that the zero map is usually
considered to be a character.  However
using
the standard terminology would require
continual use of the phrase \lq\lq non-zero
character".
\par
The set of all characters is the complement of
the zero map in the real vector space
${\rm Hom}(G,R)$.  For any character
$\chi$ let $[\chi] = \{r\chi~|~0 < r \in R\}$
be a ray in  ${\rm Hom}(G,R)$; the set
of all such rays is denoted $S(G)$.
  Since any character
of $G$ must factor through the abelianization of $G$,
 if $n$ is the integral rank of $G^{ab}$
 then $S(G) \simeq S^{n-1}$.  In particular for
a graph group $G\cG$, $S(G\cG) \simeq S^{|V(\cG)| - 1}$.
\par
Any character $\chi$ whose corresponding ray $[\chi]$
intersects an integral point of ${\rm Hom}(G,R)$ is
a {\it rational} character.  It is easy to check that the
image of a rational character is an infinite cyclic group.
\par
The Bieri-Neumann-Strebel invariants $\Sigma^k(G)$
($k$ an integer greater than zero) are open
subsets of $S(G)$ where $\Sigma^1(G) \supseteq \Sigma^2(G) \supseteq
\cdot \cdot \cdot \supseteq \Sigma^n(G) \supseteq \cdot \cdot \cdot$.
We simplify the terminology and
call these sets the {\it BNS-invariants}.
\par
The invariant $\Sigma^1(G)$ has a fairly simple description.
  Choose a finite generating set for $G$ and let
$\cC$ be the corresponding Cayley graph.  Then any character
$G ~{\buildrel \chi \over \longrightarrow}~ R$ extends to a
$G$-equivariant map
$\cC~ {\buildrel {\tilde \chi} \over \longrightarrow}~ R$
coinciding with $\chi$ on the vertices and is extended linearly
over the edges of $\cC$.  Let $\cC_+(\chi)$ be the maximal subgraph of
$\cC$ contained in ${\tilde \chi}^{-1} ([0,\infty))$.
\md
{\smps Definition.}  A character
$G~ {\buildrel \chi \over \longrightarrow}~ R$
represents a point in $\Sigma^1(G)$ if and only if
$\cC_+(\chi)$ is connected.
\md
Although it is not apparent, this definition is independent
of choice of finite generating set.
\par
The higher BNS-invariants are more difficult to define,
but we will not need the full definition for the results
in this paper.  Instead we present a characterization
which implies that a character $\chi$
 represents a point in $\Sigma^k(G)$ for $k > 1$.
\par
Assume that
 $G$ admits a $K(G,1)$ with finitely many cells in each dimension
and let $K$ be such a $K(G,1)$ which, without loss of
generality, we assume has a single vertex.
Then as in the Cayley graph example, every map $G~
{\buildrel \chi \over \longrightarrow}~ R$ extends to a $G$-equivariant
map ${\tilde K}~ {\buildrel {\tilde \chi} \over \longrightarrow}~ R$
which is defined by $\chi$ on the vertices and extended linearly
over the cells of $K$.
Let ${\tilde {K}}_+(\chi)$ be the the maximal subcomplex of ${\tilde K}$ contain
in ${\tilde \chi}^{-1}([0,\infty))$.
\md
{\smps Partial Description.}  If ${\tilde K}_+(\chi)$ is
$(n-1)$-connected, then $\chi$ represents a point
in $\Sigma^n(G)$.
\md
The higher BNS-invariants measure higher dimensional
topological properties of the kernels of maps to free
abelian groups.  Recall the following definition due
to C.T.C. Wall:
\md
{\smps Definition.}  A group $G$ has property $F_n$ if and
only if there exists a $K(G,1)$ which has finite $n$-skeleton.
Thus $F_1$ is equivalent to finite generation, $F_2$ is equivalent
to finite presentation, and a group is $F_\infty$ if it is
$F_n$ for all $n$.  For readers more versed in homological
properties, $F_n$ implies $FP_n$ for each $n$.  For background on
the $FP_n$ properties see \brown.
\md
\proclaim Theorem 1.1.  {\rm (Bieri, Neumann, Renz, Strebel)}
 Let $H$ be a normal subgroup of
$G$ with $G/H$ an abelian group of integral rank $n$.
Define $S(G,H) = \{[\chi]~|~\chi(H) = 0\} \simeq S^{n-1}$.
Then $H$ is $F_k$ if and only if $S(G,H) \subseteq
\Sigma^k(G)$.
\hfill $\square$
\par
\proclaim Corollary 1.2.  Let $\chi$ be a rational character of the
 group $G$.  Then the kernel of $\chi$ is $F_k$ if and only
if $[\chi]$ and $[-\chi]$ are both contained in $\Sigma^k(G)$.
\hfill $\square$
\par
The proof of the above theorem and corollary when $k = 1$ is
due to Bieri, Neumann and Strebel and may be found in \bns.
For the higher invariants the proof is contained in the
thesis of  Renz, and a
good discussion may be found in either \br\ or \bs.
\par
Since the finiteness properties $F_k$ are preserved when passing
to subgroups of finite index, we state most of our results in
terms of characters, not normal subgroups with abelian quotient.
In most cases
we leave the interested reader the task of translating the
theorems into the language of normal subgroups.
\vskip 1cm
\cl {\it {\rm 2.} Presentations of graph groups}
\bigskip
The presentation of a graph group given in the introduction, with
generators corresponding to the vertices of a graph $\cG$ and
relations corresponding to the edges, we call the {\it standard
presentation}.  It is intuitively the natural presentation to use;
however, it is not always the optimal presentation for specific
problems.  In \vw, the second author presented biautomatic structures
and finite complete rewriting systems for graph groups, which use an
alternate system of generators and relations.  (A similar system of
generators and relations was used in \hm\ in the more general case of
a graph {\it product} of groups.)  We give this presentation
as part of a finite complete rewriting system below.
\par
If $G$ is a group generated by a finite set $X$, adjoin a set of
formal inverses to obtain a set of monoid generators $\cB=X\cup
{\overline{X}}$ for $G$.   In general, if ${\cal S}$ is any set,
${\cal S}^*$ denotes the free monoid on ${\cal S}$.
 The inclusion of $\cB$ into $G$ extends to a
monoid homomorphism $\cB ^*~{\buildrel \mu \over \longrightarrow}~ G$.
  We will often
regard $\cB$ as a subset of $G$, suppressing the homomorphism $\mu$.
We denote the length of a word $w$ in the free monoid $\cB^*$ by $l(w)$.
\par
Let $\cC=\cC(G,X)$ denote the Cayley graph of $G$ with respect to the
generating set $X$.  Then each word $w=x_1x_2\dots x_n \in \cB^*$,
can be identified with an associated path in $\cC$,
$[0,\infty)~{\buildrel w \over \longrightarrow}~ \cC$
with $w(0)=1$, $w(i)=\mu
(a_1 ... a_i)$ for each integer $i\le n$ and $w(i)=w(n)$ for $i>n$,
where
 $[i,i+1]$ is mapped to the edge connecting $w(i)$ to $w(i+1)$.
\par
Given a generating set $X$ for $G$, a set of {\it normal forms}
for $G$ with respect to $X$ is a subset of $\cB^*$ which
bijects under $\mu$ to $G$.  In other words, a set of
normal forms defines a canonical way of representing
each group element in terms of the generators and their
inverses.
\par
A finite complete rewriting system is essentially a finite set
of rules which converts any given expression of a group
element in terms of the generators and their inverses
into the normal form for the group element.   One simply
replaces any occurance of the left-hand side of a rule in
a given word $w$ with
the right side of the rule.  If this set of rules is a finite
complete rewriting system, the above process will terminate in
a unique normal form.  The set of rules  gives
a set of relations defining the group.  Since we do not
use detailed information about finite complete rewriting systems,
the reader unfamiliar with them is
directed to \cehpt.
\par
 The normal forms for a graph group $G\cG$ used in
this paper are induced by those in \vw\  and \hm.
In order to define the normal forms we need to  define
an alternate set of generators for  a graph group.
The set
$$
\cD=\{[\{x_1,x_2,\dots ,x_n\}]\ |\ \forall i,j,\ x_i\in
X\cup{\overline{X}},\ x_ix_j=x_jx_i,\ x_i\ne {\overline{x_j}} \}.
$$
consists of a generator for each collection of standard generators
(and their inverses) which correspond to a clique in the graph $\cG$,
with the restriction that a generator and its inverse do not both occur.
If $X$ is the set of standard generators, and $\cB =X\cup
{\overline{X}}$, then we can regard $\cB \subseteq \cD$.  Further there is a
map $\zeta$ such that $\cD^*~{\buildrel\zeta\over\longrightarrow}~
\cB^* ~{\buildrel\mu\over\longrightarrow}~ G\cG$, where $\zeta$ is
determined by a total ordering of $X$: if $x_{i+1}<x_i$, then
$\zeta([\{x_1,x_2,\dots ,x_n\}])=x_nx_{n-1}\dots x_1$.  Since $\cB$ is
contained in the image of $\cD$ under $\zeta$, the composition $\mu
\circ \zeta$ maps $\cD^*$ onto $G\cG$.
\par
Let
$w_1,w_2\in\cD^*$, $x\in X\cup {\overline{X}}$, and
$+$ denote disjoint union.  The rewriting rules for elements of
$\cD^*$ are:
\md
 1. $w_1[u][v+\{x\}]w_2 \rightarrow
w_1[u+\{x\}][v]w_2$
\ \ \ provided $\forall
y\in u$, $x$ and $y$ commute, $x\ne y$, and $x\ne {\overline{y}}$.
\par
2. $w_1[u+ \{x\}][v+
\{{\overline{x}}\}]w_2\rightarrow w_1[u][v]w_2$.
\par
 3. $w_1[\emptyset]w_2\rightarrow w_1w_2$.
\md
This is a ``left greedy'' reduction system; we simply move
an element of $\cB$ to the left whenever possible, cancelling if
necessary (as in (2)).  That this is a finite complete
rewriting system for the graph group is shown in \vw.
The presentation of $G\cG$ with generating set $\cD$
and relations given by the rewriting system is the
{\it clique presentation}.
\md
We prefer to work with the standard presentation.  However, the
rewriting system as defined only gives normal forms
for group elements in terms of the clique generators.
Using the map $\zeta$ these normal forms can
be easily converted to normal forms in $\cB^*$.
If $w$ is any word in $\cB^*$, then
$w$ can be viewed also as a word in $\cD^*$ consisting of singletons.
Applying the rewriting rules to $w$ then gives
 an element $[u_1][u_2]\dots[u_k]\in \cD^*$ in normal form.
This
induces a normal form  $\zeta ([u_1])\zeta ([u_2])\dots\zeta
([u_n])\in\cB^*$.  It is this  set of normal forms which
will be used.
\par
If $w$ and $w'$ are elements of $\cB^*$ we use
$w \longrightarrow^* w'$ to indicate that one can get
from $w$ to $w'$ by applying rewriting rules to $w$
(thought of as a word in $\cD^*$) and then using the
map $\zeta$ to convert the rewritten word into a word
in $\cB^*$.
\vskip 1cm
\cl {\it {\rm 3.} Normal forms and characters}
\bigskip
In determining the structure of $\Sigma^k(G\cG)$
it suffices to discuss only those characters which
map the standard generators to the intersection of
$S(G\cG)$ with the cone formed by the
positive coordinate axes.  Specifically  we have the
following result.
\par
\proclaim Proposition 3.1.  Let $\Sigma_P^k(G\cG)
\subseteq S(G\cG) \simeq S^{n-1}$ be
the intersection of $\Sigma^k(G\cG)$ with the closed
 cone formed by the positive coordinate axes.
   Let $C$ be the finite Coxeter group
generated by  reflections in the $n$
 hyperplanes determined
by each linearly indepedent set of
$n-1$ coordinate axes.  Then the image of $\Sigma_P^k(G\cG)$
under the action of $C$ is $\Sigma^k(G\cG)$.
\par
{\it Proof.} Choose any point in $S(G\cG)$ and let
$G\cG~{\buildrel \chi \over \longrightarrow}~ R$ be a character
representing the point.  Using the standard presentation of $G\cG$,
$\chi$ is determined by the images of the group elements corresponding
to the vertices of $\cG$.  Let $ G\cG~ {\buildrel \chi' \over
\longrightarrow}~ R$ be the map which takes $v$ to $-\chi(v)$ for each
generator $v$ with $\chi (v)<0$, and which agrees with $\chi$ on all
other generators.
\par
Changing from $\chi$ to $\chi '$ is essentially the same as replacing
$\{v\ |\ \chi(v)<0\}$ in the generating set by $\{{\overline{v}}\ |\
\chi(v)<0\}$.  But changing a generator from $v$ to ${\overline {v}}$,
and similarly changing all the appearances of $v$ in the relations,
gives an isomorphic presentation of $G\cG$.  Because of this symmetry
in the defining relations, the sets ${\tilde K}_+(\chi)$ and
${\tilde K}_+(\chi')$ are
isomorphic, so $[\chi ']\in\Sigma^k(G\cG)$ if and only if $[\chi] \in
\Sigma^k(G\cG)$.  Since each of the generating reflections of $C$ simply
maps a single generator $v$ to ${\overline{v}}$, there exists $c\in C$
such that $c[\chi ']=[\chi ]$.  Hence $\Sigma^k(G\cG) = C\cdot
\Sigma_P^k(G\cG)$.
\hfill $\square$
\md
{\smps Definition.}
If $\chi$ is a character of a graph group $G\cG$
then $\cL(\chi)$ is the full subgraph of $\cG$
generated by the vertices corresponding to generators
which are {\it not}
mapped to $0$ by $\chi$.  We refer to this
 as the \lq\lq living" subgraph of
$\cG$.
In view of Proposition 3.1., we may assume $\chi(x)>0$ for each generator
$x$ corresponding to a vertex in $\cL(\chi)$.
\md
In section 4 we discuss the connectedness of $\cC_+(\chi)$.
We use the set of normal forms described in section 2 to
give some control over the form of paths in $\cC_+(\chi)$.
In particular we will need to refer to the following two lemmas.
\md
\proclaim Lemma 3.2.  Let $\chi$ be a character for a graph
group $G\cG$.  Let $w$ be a path
in $\cC_+(\chi)$ of minimal length from 1 to $\mu(w)$ and
assume $\mu(w)$ has normal form $z$.  If $x$ is any generator which
occurs in $w$ but not in $z$, then $x$ corresponds to
a vertex in $\cL(\chi)$ and we can write $w = pxq{\overline {x}}r$.
\par
 {\it Proof.} Choose a sequence of reductions $w \rightarrow ^* z$.
Since $x$ doesn't occur in $z$, either $w=pxq{\overline {x}}r$ or
$w=p{\overline {x}}qxr$, where these ``visible'' $x$'s cancel
somewhere in this sequence.
Let $w'=pqr$.  Then $w' \rightarrow ^* z$ via the same sequence of
reductions except for the steps involving these $x$'s.  So $\mu
(w)=\mu (w')$.
\par
Suppose $x$ is a generator corresponding to a vertex $v \not\in
\cL(\chi)$.  Then for any decomposition of $q$ as $q=q'q''$,
$\chi (pq') = \chi
(pxq') = \chi (p{\overline{x}}q')\ge 0$.  Similarly, if $r=r'r''$,
then $\chi (pqr') = \chi (pxq{\overline{x}}r') = \chi
(p{\overline{x}}qxr')\ge 0$.  Thus every initial segment of $w'$ is
contained in $\cC_+(\chi)$, and hence the word
 $w'$ is a path in $\cC_+(\chi)$ from
1 to $\mu (w)$.  But $w'$ is shorter than $w$
 which contradicts the minimality of $w$.  Hence,
$x$ corresponds to a vertex in $\cL(\chi)$.
\par
Now suppose $w=p{\overline{x}}qxr$ with $\chi(x)>0$.  If
$q=q'q''$, then $\chi
(pq')>\chi (p{\overline{x}}q')\ge 0$ and if $r=r'r''$, then $\chi
(pqr')=\chi (p{\overline{x}}qxr')\ge 0$.  So again $w'$ is a path in
$\cC_+(\chi)$ from 1 to $\mu (w)$, a contradiction.
\hfill $\square$
\md
\proclaim Lemma 3.3.  Let $w$ be as in the previous lemma.  Assume
$w \rightarrow ^* w' \rightarrow ^* z$ where $z$ is the normal form for
$\mu(w)$
and $l(w') > l(z)$.  Then there exits
a generator $x$, corresponding to a vertex in
$\cL(\chi)$, such that
$w' = p'xq'{\overline {x}}r'$ where $x$ commutes with all the
letters in $q'$.
\par
 {\it Proof.} Since $l(w')>l(z)$ and $w'\rightarrow^*z$,
a length-reducing rule must occur.  Hence, there exists
$x\in X$ and words $p'$, $q'$, $r'$, with either
$w'=p'xq'{\overline{x}}r'$ or $w'=p'{\overline{x}}q'xr'$, where $x$
commutes with all the letters in $q'$.  It remains to show that
$\chi(x) > 0$ and that $x$ occurs before ${\overline {x}}$.
\par
Since $w \rightarrow^* w'$, there exist words
$p$, $q$, and $r$ with either
$w=pxq{\overline{x}}r$ or $w=p{\overline{x}}qxr$ ($x$ does not
necessarily commute with all the letters in $q$).  As in the proof of
Lemma 3.2, if $x$ does not correspond to a vertex in $\cL(\chi)$ or if
$w=p{\overline{x}}qxr$, then $pqr$ defines a path in $\cC_+(\chi)$ from 1 to
$\mu (w)$, a contradiction.
\hfill $\square$
\vskip 1cm
\cl {\it {\rm 4.} $\Sigma^1(G\cG)$}
\bigskip
Whether or not a character $\chi$ of a graph
group represents a point in the BNS-invariant
$\Sigma^1(G\cG)$ is determined by the living subgraph
$\cL(\chi)$ of $\cG$.  Recall that a subgraph $\cL$ of
$\cG$ is
{\it dominating} if for every $v \in \cG\backslash\cL$,
$d(v,\cL) = 1$.
\par
\proclaim Theorem 4.1.  Let $\chi$ be a character of
a graph group $G\cG$. The subgraph $\cC_+(\chi)$ of $\cC$
is connected if and only if $\cL(\chi)$
 is a connected dominating subgraph of $\cG$.
\par
 {\it Proof.} $(\Rightarrow$):  Assume
$A=\{a_1,a_2,\dots,a_n\}$ is the subset of the generating set
$X$,
corresponding to the vertices in $\cL(\chi)$.  Assume also that
there is a chosen total ordering on $X$ where
 $a_i<a_{i+1}$ and for convenience, $\chi(a_1) \le \chi(a_i)$ for
all $i$.  We first show that
$\cL(\chi)$ is a dominating subgraph of $\cG$.
Let $b$ be a generator with
$\chi(b)=0$, let $a_k\in A$, and let
$w$ be a path of minimal length in $\cC_+(\chi)$ from 1 to
$\mu({\overline{a_k}}ba_k)$.   If $l(w)=1$, then $a_k$ and $b$ commute,
and hence $d(b,\cL(\chi))=1$.  So assume $l(w)>1$.
\par
Since $a_k$ and $b$ don't commute, $w$ isn't in normal form (the
normal form for $w$ is ${\overline{a_k}}ba_k$, which is not in
$\cC_+(\chi)$ as a path).  By repeated use of Lemma 3.3,
we can define a sequence of words
$w=w_1,w_2,w_3,\dots, w_n={\overline{a_k}}ba_k$, where for each $j$,
$w_j=
p_ja_{i_j}q_j{\overline{a_{i_j}}}r_j$ and the rewriting
$w_j\rightarrow^* w_{j+1}$ consists of commuting letters within
$p_j, q_j$ and $r_j$ along with cancelling the $a_{i_j}$ and ${\overline
{a_{i_j}}}$ terms.  Further, it can be assumed that $a_{i_j}$ commutes
with all the letters in $q_j$.
\par
By hypothesis $\chi(b)=0$, so there is only
one occurance of $b$ in $w$ by minimality.
Suppose $b$ occurs in none of the $q_j$.  Then we can write $w=ubv$, where
${\overline{a_k}}$ occurs in $u$ and all the other letters in $u$
occur in inverse pairs.  Thus $\chi(u)=-\chi(a_k)<0$, which
contradicts the assumption that $w$ is a path in $\cC_+(\chi)$.
\par
It follows that there is a $j$ such that $b$ occurs in $q_j$, and therefore $b$
commutes with $a_{i_j}$.  Since the distance from $b$
to $\cL(\chi)$ is $1$,
 $\cL(\chi)$ is a dominating
subgraph of $\cG$.
\md
It remains to show that $\cL(\chi)$ is connected.
 For each $i\ne 1$, let $w(i)$ be a path of
minimal length in $\cC_+(\chi)$ from 1 to $\mu({\overline{a_1}}a_i)$.
The vertex $\mu({\overline{a_1}}a_i)$ is contained in $\cC_+(\chi)$
since we have assumed that $\chi(a_1) \le \chi(a_i)$ for
all $i$.
Note that the normal form for $w(i)$ is ${\overline{a_1}}a_i$, which
doesn't lie in $\cC_+(\chi)$ as a path since we may
assume that $\chi$ maps the generators in $A$ to positive
real numbers.  The proof of the following
claim then concludes the proof of $(\Rightarrow)$.
\md
{\smps Claim:} For each $i$, there exists a path in $\cG$ from $a_1$ to
$a_i$.
\md
{\it Proof of Claim.}  We proceed by induction on the length of
$w(i)$.
\md
\underbar{Base}:  $l(w(i))=2$.  Then we must have
$w(i)=a_i{\overline{a_1}}$, so $a_1$ and $a_i$ commute.  Hence $a_1$
and $a_i$ are adjacent in $\cG$.
\md
\underbar{Inductive Step}:  $l(w(i))>2$.  In this case $a_1$ and
$a_i$ do not commute.  Write $w(i)=ua_iv$, where this is the first
occurance of $a_i$.  Then $v$ must be non-empty since otherwise
$\chi(u)=-\chi(a_1)<0$, which contradicts the assumption that $w$ lies
in $\cC_+(\chi)$.
\par
By Lemma 3.2, if there is an $a_k$ which
occurs in $v$, then so does ${\overline{a_k}}$,
and $v=pa_kq{\overline{a_k}}r$.  Suppose such is the case, and that
$a_k$ is the last element of $A$ in $v$, i.e., $q,r\in ({\overline{A}})^*$.
Let $w'=ua_ipqr$ and notice that $\mu(w')=\mu(w)$.
As in the proofs of Lemmas 3.2 and 3.3 it can be shown
that $w'$ is also contained in $\cC_+(\chi)$ and is of
shorter length than $w$, contradicting the minimality
of $w$.
\par
We may then assume that $w(i)=ua_iv$ where $v\in ({\overline{A}})^*$.
  It follows that each
element of ${\overline{A}}$ which occurs in $v$ commutes with $a_i$.  Let
$v=v_1{\overline{a_j}}$.  Then
$\mu(ua_iv_1{\overline{a_j}})=\mu(uv_1{\overline{a_j}}a_i)$ and
$\mu(uv_1)=\mu({\overline{a_1}}a_j)$.
Also, $\chi(uv_1)=\chi(a_j)-\chi(a_1)\ge 0$.
Thus for each decomposition $v_1=v_1'v_1''$,
$\chi(uv_1')=\chi(uv_1)+\chi({\overline{v_1''}}) \ge \chi(uv_1)\ge 0$.
\par
Hence $uv_1$ is a path in $\cC_+(\chi)$ from 1 to $\mu (uv_1)=\mu
({\overline{a_1}}a_j)$ of length $l(w(i))-2$.  By definition,
$l(w(j))\le l(uv_1)$, so by the inductive hypothesis there exists a
path in $\cG$ from $a_1$ to $a_j$.  Since $a_i$ and $a_j$ commute
(${\overline{a_j}}$ occured in $v$), we
have a path from $a_1$ to $a_i$ in $\cG$.
\hfill $\square$
\md
 $(\Leftarrow$):  Let $g$ be a vertex in $\cC_+(\chi)$; we will
construct a path in $\cC_+(\chi)$ from 1 to $g$.  First, we need some
notation.  We will be blatantly indentifying the set of standard
generators with $V(\cG)$ and at times we will use the symbol
\lq\lq $\chi$" to denote the  map induced on the free monoid
$\cB^*$ by the character $\chi$ in ${\rm Hom}(G\cG, R)$.  By Proposition
3.1, we can assume $\chi (x)>0$ for all standard generators $x$.
\par
Let $y_1,y_2\in V(\cG)\cup {\overline{V(\cG)}}$.  Since $\cL(\chi)$ is
connected and dominating, there exists a path $y_1,x_1,x_2,\dots
,x_k,y_2$ from $y_1$ to $y_2$ in $\cG$ with $\chi (x_i)>0$ for all
$i$.  For positive integers $r$ and $s$, define
$$
p(y_1^r,y_2^s)=x_1^{n_1}y_1^rx_2^{n_2}{\overline{x}}_1^{n_1}
x_3^{n_3}{\overline{x}}_2^{n_2}\dots
x_{k-1}^{n_{k-1}}{\overline{x}}_{k-2}^{n_{k-2}}x_k^{n_k}
{\overline{x}}_{k-1}^{n_{k-1}}y_2^s
$$
where the positive integers $n_i$ are minimal such that
$\chi(x_1^{n_1})\ge |\chi(y_1^r)|$,
$\chi(x_i^{n_i})\ge\chi(x_{i-1}^{n_{i-1}})$ for $2\le i\le n-1$, and
$\chi(x_k^{n_k})\ge\chi(x_{k-1}^{n_{k-1}})+|\chi(y_2^s)|$.  The net
effect of these choices of $n_i$ is that the path defined by
$p(y_1^r, y_2^s)$ starting at $1$ in the Cayley graph is contained
in $\cC_+(\chi)$.  In addition,
define $z(y_1^r,y_2^s)={\overline{x}}_k^{n_k}$.  The path $p(y_1^r, y_2^s)$
is chosen so that much collapsing occurs, in particular
 $\mu (y_1^ry_2^s)=\mu(p(y_1^r,y_2^s)z(y_1^r,y_2^s))$.
\par
Now, let $w_0\in (V(\cG)\cup {\overline{V(\cG)}})^*$.
Then $w_0=uv_1v_2\dots v_m$ where $v_i\in V(\cG)\cup
{\overline{V(\cG)}}$, $u\subseteq\cC_+(\chi)$, and
$\mu(uv_1)\not\in\cC_+(\chi)$.  (If $w=u$, then we are done.)  Using
the machinery in the previous paragraph, replace $v_1v_2$ by
$p(v_1,v_2)z(v_1,v_2)=p_1z_1$, forming $w_1=up_1z_1v_3v_4\dots v_m$.
Then $\mu(w_1)=\mu(w_0)$ and $up_1\subseteq\cC_+(\chi)$.  Similarly, replace
$z_1v_3$ by $p(z_1,v_3)z(z_1,v_3)=p_2z_2$, forming
$w_2=up_1p_2z_2v_4\dots v_m$, so that $\mu(w_2)=\mu(w_1)=\mu(w_0)$ and
$up_1p_2\subseteq \cC_+(\chi)$.  Continuing in this manner, we arrive
at the word $w_m=up_1p_2\dots p_{m}z_{m}$, where $\mu(w_m)=\mu(w_0)$ and
$up_1p_2\dots p_{m}\subseteq\cC_+(\chi)$.
\par
By construction, $z_{m}=y^{\alpha}$ for some $y\in V(\cG)\cup
{\overline{V(\cG)}}$ and some non-negative integer $\alpha$.
Since both $\mu(w_m)$ and $\mu (up_1p_2\dots p_{m})$
are vertices in $\cC_+(\chi)$, if
$\chi (y)\ge 0$, then $\chi(up_1p_2\dots p_{m}y^i)\ge 0$ for $1\le
i\le\alpha$, and if $\chi(y)<0$, then $\chi(up_1p_2\dots
p_{m}y^i)=\chi (w_m{\overline{y}}^{\alpha-i})\ge 0$ for $1\le
i\le\alpha$.   It follows that $w_m$ is a path in
$\cC_+(\chi)$ from 1 to $g$, and hence $\cC_+(\chi)$ is connected.
\hfill $\square$
\md
{\smps Example 1:}  The free group on two generators
can be thought of as the graph group corresponding to the
null graph with two vertices.  If $\chi$ maps each generator
to $1 \in Z$, then $\cL(\chi)$ is disconnected, hence
by Theorem 4.1 the subgraph $\cC_+(\chi)$ is disconnected.
\par
In the figure below, the gray scale indicates the parts of the
Cayley graph which are in $\cC_+(\chi)$, and it is apparent
that $\cC_+(\chi)$ is disconnected.
\vfill
\cl {\smps Fig. 1}
\eject
{\smps Example 2:} If we add an additional vertex and two edges to the
graph in example 1, forming a graph $\cG = x$---$y$---$z$, then $G\cG$
is the direct product of an infinite cyclic group with a free group.  Let
$\chi$ be the character sending each generator to $1 \in Z$.  (Thus
$\chi$ restricted to the subgroup generated by $x$ and $z$ is the same
map as was discussed in example 1.)
\par
By Theorem 4.1, $\cC_+(\chi)$ should be connected.  Below
is part of the Cayley graph, with the $x$ and $z$ generators
forming the tree in front and the positive $y$-direction going
back into the page.  Once again the gray scale indicates
the subgraph $\cC_+(\chi)$.
\vskip 8cm
\cl {\smps Fig. 2}
\vskip 1cm
\cl {\it {\rm 5.}  The simplicial structure of $\Sigma^1(G\cG)$}
\bigskip
 The unit sphere $S^{n-1}$ admits a natural simplicial
structure where the vertices correspond to the intersection
of the axes with $S^{n-1}$, two vertices are joined
by an edge if they are a distance $\ds {\pi \over 2}$
apart and any complete graph with $m$ vertices is
filled by an $(m-1)$-simplex.  Equivalently, it is the
complex formed by stating that any $m$
vertices corresponding to linearly independent
(positive or negative) coordinate axes are
the vertices of  an $(m-1)$-simplex.
\md
{\smps Notation.}
Let $\cS(G\cG)$ denote the sphere
 $S(G\cG)$ with the simplicial structure
sketched above.
\md
As was mentioned in section 1,
for any finitely presented group $G$ the set
$\Sigma^1(G)$ is always an open subset of $S(G)$.
By the previous theorem it follows that $\Sigma^1(G\cG)$
actually is the union of open simplices
 in the simplicial
decomposition $\cS(G\cG)$ (see Theorem 5.1 below).
The vertices of $\cS(G\cG)$ correspond to characters
of the graph group $G\cG$ whose living subgraph
consists of a single vertex.  Thus a vertex
of $\cS(G\cG)$ is included
in $\Sigma^1(G\cG)$ if and only if the corresponding
vertex in $\cG$ dominates $\cG$.  Similarly the
 open edges of $\cS(G\cG)$ correspond to pairs of vertices
in $G\cG$, and the entire open edge is included
in $\Sigma^1(G\cG)$ if and
only if the  living subgraph generated by the
pair of vertices is connected
and dominating.   The conditions for the higher
dimensional open simplices are similar.
\md
{\smps Definition.}  Since any open simplex  of
$\cS(G\cG)$ is defined by a collection of coordinate
axes, each open simplex defines a collection
of vertices of $\cG$.
In an abuse of terminology, if $\sigma$ is an open
simplex of $\cS(G\cG)$ then
 the {\it living subgraph}  $\cL(\sigma)$ is
the full subgraph of $\cG$ generated by the vertices
of $\cG$ which correspond to the coordinate axes
defining $\sigma$.
\md
\proclaim Theorem 5.1.  For $G\cG$ a graph group,
the BNS-invariant $\Sigma^1(G\cG)$
is a union of open simplices in $\cS(G\cG)$.
In particular, an open simplex $\sigma$ of $\cS(G\cG)$
is contained in $\Sigma^1(G\cG)$ if and only if the
corresponding living subgraph
$\cL(\sigma)$ is connected and dominating.
\par
 {\it Proof.}   If $\chi$ and $\chi'$ are two characters defining
points in the same simplex $\sigma$ of $\cS(G\cG)$, then
$\cL(\chi) = \cL(\chi') = \cL(\sigma)$.  Thus all the points
in $\sigma$ have the same living subgraph of $\cG$, hence either
the entire open simplex $\sigma$ is contained in
$\Sigma^1(G\cG)$ or $\sigma \cap \Sigma^1(G\cG) =
\emptyset$.
\hfill $\square$
\md
\proclaim Corollary 5.2.  For every graph group $G\cG$,
$\Sigma^1(G\cG)$
is rational polyhedral.
If $\cG$ is disconnected then $\Sigma^1(G\cG)$
is empty, and if $\cG$ is connected then
the closure of $\Sigma^1(G\cG)$ is $S(G\cG)$.
\par
 {\it Proof.}  That $\Sigma^1(G\cG)$ is rational polyhedral~---~that
is, can be defined by finitely many inequalities with integer
 coefficients~---~follows directly from Theorem 5.1 and the description of
$\cS(G\cG)$.
\par
 If $\cG$ is disconnected then there
exist no connected dominating subgraphs, hence
no simplex of $\cS(G\cG)$ can be contained in
$\Sigma^1(G\cG)$.
\par
If $\cG$ is connected, then the graph $\cG$ is
a connected dominating subgraph, so
every maximal simplex in $\cS(G\cG)$ will
be contained in $\Sigma^1(G\cG)$.  The result
follows since the closure of the maximal simplices
of $\cS(G\cG)$ is $S(G\cG)$.
\hfill $\square$
\md
It is an open question whether the
invariant $\Sigma^1(G)$ is rational
polyhedral for all finitely generated groups $G$.
\par
\proclaim Proposition 5.3.  Let $\tau$ be an open simplex
of $\cS(G\cG)$ whose closure contains an open simplex
$\sigma$ which is contained in $\Sigma^1(G\cG)$.  Then
$\tau$ is contained in $\Sigma^1(G\cG)$.
\par
 {\it Proof.}  The living subgraph $\cL(\sigma)$ is
connected and
dominating by hypothesis.  Since $\cL(\tau)$
 only adds vertices to the graph
$\cL(\sigma)$ it is clear that $\cL(\tau)$
is also dominating.  Further, since $\cL(\sigma)$
is dominating, every  vertex  $v \in \cL(\tau) \backslash \cL (\sigma)$
is a distance $1$ from $\cL(\sigma)$.  Hence the
full subgraph generated by adding these vertices
to $\cL(\sigma)$ is also connected.
\hfill $\square$
\md
\proclaim Corollary 5.4.  For a graph group $G\cG$, $\Sigma^1(G\cG) \simeq
S^{|V(\cG)| - 1}$ if and only if $\cG$ is a complete graph.
\par
 {\it Proof.}  By the previous proposition, $\Sigma^1(G\cG)
\simeq S^{|V(\cG)|-1}$ if and only if the vertices of $\cS(G\cG)$
are contained in $\Sigma^1(G\cG)$.
The vertices  are contained
in $\Sigma^1(G\cG)$ if and only if
each corresponding vertex in $\cG$ is
 dominating.
The corollary follows since each vertex
dominating is equivalent to the graph being a
complete graph.
\hfill $\square$
\md
{\smps Example:} If $\cG$ is given by $x$---$y$---$z$,
then $\Sigma^1(G\cG)$ is a 2-sphere less the great circle determined by
$x$ and $z$.  It follows from Theorem 8.3
 that $\Sigma^k(G\cG)$ is a 2-sphere with a great circle removed
for all $k$.
\vskip 1cm
\cl {\it {\rm 6}. Free abelian quotients}
\bigskip
Because of the symmetry discussed in Proposition 3.1, $[\chi] \in
\Sigma^1(G\cG)$ if and only if $[-\chi] \in \Sigma^1(G\cG)$.  Thus
Theorem 6.1 follows from Corollary 1.2 and Theorem 4.1.
\par
\proclaim Theorem 6.1.  Let $\chi$ be a rational character
of a graph group $G\cG$.  The kernel of $\chi$ is
finitely generated if and only if $\cL(\chi)$ is a
connected and dominating subgraph of $\cG$.
\hfill $\square$
\par
If $G\cG~ {\buildrel \chi \over \longrightarrow}~ Z^n$ with $n>1$ then
the situation is more complex.  Let $H$ be the kernel of $\chi$.  Then
$S(G\cG,H)$ is naturally identified with an $(n-1)$-sphere as was
mentioned in Theorem 1.1.  The sphere arises since any non-zero map to
a free abelian group can be extended to a collection of characters via
the elements of ${\rm Hom}(Z^n,R)$ and the non-zero elements of ${\rm
Hom}(Z^n,R)$, scaled by the positive reals, is naturally identified
with a copy of $S^{n-1}$.
\par
There is the following partial result.  Recall that a graph $\cG$ is
$m$-{\it connected\/} if  there are $m$ vertices whose
removal yields a disconnected graph.
\md
\proclaim
Proposition 6.2.  If the graph $\cG$ is $m$-connected, then $G\cG$
admits no map onto $Z^n$ for $n>m$ with finitely generated kernel.
\par
 {\it Proof.}  Let $\chi$ map $G\cG$
onto $Z^n$ for some $n>m$.  Let $v_1,v_2,\dots ,v_m$ be the vertices
whose removal disconnects $\cG$, and let $x_1,x_2,\dots ,x_m$ be the
corresponding generators of $G\cG$.  Since $\chi$ is surjective, we
can choose a non-zero element $\phi\in{\rm Hom}(Z^n,Z)$ such that
$(\phi\circ\chi )(x_i) = 0$ for $i=1,2,\dots ,m$.  Since each
$v_i\not\in\cL(\phi\circ\chi)$, $\cL(\phi\circ\chi)$ is not connected,
so $[\phi\circ\chi]\notin\Sigma^1(G\cG)$.  Clearly $[\phi\circ\chi]\in
S(G\cG,Ker(\chi))$, so $S(G\cG,Ker(\chi))\not\subseteq\Sigma^1(G\cG)$
and thus $Ker(\chi)$ is not finitely generated by Theorem 1.1.
\hfill $\square$
\md
By Theorem 5.1 it is clear that an open simplex of $\cS(G\cG)$ is
contained in $\Sigma^1(G\cG)$ if and only if some rational point in
the simplex is contained in $\Sigma^1(G\cG)$.  Thus it suffices to
consider characters defined by composing $\chi$ with maps in ${\rm
Hom}(Z^n,Z)$.
\par
\par
\proclaim Proposition 6.3.  A map $G\cG ~{\buildrel \chi \over
\longrightarrow}~ Z^n$ has a finitely generated kernel
if and only if for every nonzero $\phi \in {\rm Hom}(Z^n,Z)$,
$\cL(\phi\circ\chi )$ is connected
and dominating.
\par
{\it Proof.}  If there exists nonzero $\phi$ with $\cL (\phi\circ\chi )$
non-connected or non-dominating, then it follows as in the proof of
Proposition 6.2 that $Ker(\chi)$ is not finitely generated.
Conversely, every representative of a
rational map $\psi$ in $S(G\cG,Ker(\chi))$
factors as $\psi =\phi\circ\chi$ for some $\phi$.  Since
$[\phi\circ\chi ]\in\Sigma^1(G\cG)$ by hypothesis,
$S(G\cG,Ker(\chi))\subseteq\Sigma^1(G\cG)$, whence
$Ker(\chi)$ is
finitely generated by Theorem 1.1.
\hfill $\square$
\md
The collection ${\rm Hom}(Z^n,Z)$ is obviously a
large set of maps,  even modulo scalar multiplication.
It is not necessary, however, to
check infinitely many maps to establish that the kernel
of a map onto a free abelian group is finitely generated.
\md
{\smps Definition.}  Given a map $G\cG~ {\buildrel \chi \over
\longrightarrow}~ Z^n$
from a graph group to a free abelian group, we call a set
of standard generators $\{ x_1, \dots , x_n\}$ {\it $\chi$-linearly
independent} if for each subset $I$ of $\{1,2 \dots ,n\}$
and each $j \in \{1,2 \dots ,n\}\backslash I$,
the abelian group $\langle~\chi(x_i)~|~
 i \in I~\rangle$ is of infinite index in
$\langle~\chi(x_i)~|~i \in I \cup \{j\}~\rangle$.
A map $\phi \in {\rm
Hom}(Z^n,Z)$ is $\chi$-{\it oriented\/} provided
 $n-1$ $\chi$-linearly independent
standard generators of $G\cG$ are contained in
the kernel of $\phi \circ \chi$.
\par
\proclaim Theorem 6.4.  A map $G\cG~ {\buildrel \chi \over
\longrightarrow}~ Z^n$ has a finitely generated kernel if and only
if for every $\chi$-oriented map $\phi \in {\rm Hom}(Z^n,Z)$,
$\cL(\phi\circ\chi)$ is connected and dominating.
\par
 {\it Proof.}  If the kernel is finitely generated, then the result
follows from 6.3.  Conversely, if $\phi \in {\rm Hom}(Z^n,Z)$, then
there exists a $\chi$-oriented map $\phi'$ with
$\cL(\phi'\circ\chi)\subseteq\cL(\phi\circ\chi)$.  Since
$\cL(\phi'\circ\chi)$ is connected and dominating, so is
$\cL(\phi\circ\chi)$, and again the result follows from 6.3.
\hfill $\square$
\md
 Theorem 6.4 yields some interesting, and even easy
to apply, corollaries.
\par
\proclaim Corollary 6.5.  Let $G\cG~ {\buildrel \chi \over
\longrightarrow}~ Z^n$ and assume $V(\cG)$ is
$\chi$-linearly
independent.  Then $Ker(\chi)$ is finitely generated if and only
if $\cG$ is a complete graph.  In particular, the commutator
subgroup of a graph group is finitely generated if and only
if $\cG$ is a complete graph.
\par
 {\it Proof.}  If $\cG$ is a complete
graph, then $\cL(\phi\circ\chi)$ is connected and dominating for all
nonzero $\phi:Z^n\longrightarrow Z$.  Conversely, since $\chi(V(\cG))$
is linearly independent, for each $v\in V(\cG)$ there exists
$\phi_v:Z^n\longrightarrow Z$ which projects onto the maximal cyclic
subgroup containing $\chi (v)$, so that $\cL (\phi_v\circ\chi)=\{v\}$.
Since $Ker(\chi)$ is finitely generated, $\cL (\phi_v\circ\chi)$ must
dominate $\cG$ for each $v\in V(\cG)$ and therefore $\cG$
must be complete.  The last sentence is immediate since under
the map defined by abelianization $V(\cG)$ is a linearly
independent set.
\hfill $\square$
\md
\proclaim Corollary 6.6.  If $\cG$ admits
$n$ disjoint connected dominating sets,
then there is a map from $G\cG$
onto $Z^n$ with finitely generated
kernel.
\par
 {\it Proof.}  If $A_1,A_2,\dots,A_n$ are the connected dominating
sets, define $G\cG~{\buildrel \chi \over \longrightarrow}~ Z^n$  by sending each
vertex in $A_i$ to the standard unit vector $e_i$.  Clearly for any
$\phi:Z^n\longrightarrow Z$, $\cL(\phi\circ\chi)$ is connected
and dominating.
\hfill $\square$
\md
{\smps Example.}  The converse of
Corollary 6.6 is false.  For instance, let $\cG$
be the graph with five vertices $v_1, \dots , v_5$
and an edge joining $v_i$ to $v_{i+1}$  for each
$i$ (indices taken
modulo 5).  Let $x_i$ denote the standard generator
corresponding to the vertex $v_i$.  Using Theorem 6.4
it can be checked that the map $\chi$
from $G\cG$ onto the free abelian group of rank 2
defined by $\chi(x_1)=\chi(x_2)=(1,0)$, $\chi(x_3)=\chi(x_4) = (0,1)$
and $\chi(x_5)=(1,1)$ has a finitely generated kernel.
However $\cG$ does not contain two disjoint connected
dominating subgraphs.
\vskip 1cm
\cl {\it {\rm 7.}  Eilenberg-MacLane complexes}
\bigskip
To extend the discussion of the topological
properties of kernels of maps to free abelian groups,
we need to construct $K(G\cG,1)$ complexes
in order to apply Theorem 1.1.  The
universal covers of the complexes
we construct can  be thought of  as piecewise
Euclidean cubical complexes, that is, as complexes where each
cell is given the metric structure of a Euclidean
$n$-cube of unit side length.
\par
Let $C$ be such a cubical complex.   If $v$ is any vertex of
$C$, then the sphere of radius $1$
about $v$ inherits a natural simplicial structure.
Call this sphere the {\it link} of $v$.
\vfill
\cl {\smps Fig. 3}
\eject
{\smps Definition.} A simplicial complex $\cX$ is a
{\it flag complex} if given any complete graph
$\cC$ contained in $\cX^{(1)}$, there is some simplex
$\cal S$ in $\cX$ with ${\cal S}^{(1)} = \cC$.
\md
The proof of the following theorem is based
on the piecewise Euclidean metric structure
for cubical complexes.  If the links of all
the vertices are flag complexes, then the
cubical complex is ``nonpositively curved",
and hence is
a unique geodesic space.  For details see
section 4.2 of \gromov.
\par
\proclaim Theorem 7.1.  {\rm (Gromov)} Let $C$ be a
1-connected
cubical complex.  If the link of
each vertex is a flag complex, then $C$ is contractible.
\hfill $\square$
\par
Let $K \cG^{(2)}$ be the standard
$2$-complex for the standard presentation of $G\cG$.
The universal cover of $K \cG^{(2)}$
  is a cubical $2$-complex
which we denote $\tKG^{(2)}$.  Let $Q$
be a set of $n$ vertices
in $\cG$ which are the vertices of a
 clique in $\cG$.  Let $KQ$ be
the corresponding set of edges in $K \cG^{(2)}$.
Then $KQ$ is simply the quotient of the
1-skeleton of a Euclidean $n$-cube after
identifying all opposite faces.
In other words, $KQ$ is the $1$-skeleton
of the $n$-torus.
\par
Let $T^n$ be the $n$-torus constructed
as above by identifying
opposite faces of the
Euclidean $n$-cube.  Thus $T^n$ is a subcomplex
of $T^m$ if $n < m$.
For each  clique $Q_i$ in $\cG$
Let $T_i$ be a copy of $T^n$ where $n
= |V(Q_i)|$.   If two  cliques
$Q_i$ and $Q_j$ intersect in a (possibly trivial)
sub-clique with $k$ vertices, then identify
the corresponding copies of $T^k$ in
$T_i$ and $T_j$.  In particular,
 if the cliques $Q_i$
and $Q_j$ do not intersect, $T_i$ and $T_j$ are
joined only at the single vertex.
Call the resulting cell complex $K\cG$.
\par
\proclaim Theorem 7.2.  The cell complex $K\cG$ is
a $K(G\cG,1)$.
\par
 {\it Proof.}  The $2$-skeleton of $K\cG$ is
simply the original standard $2$-complex for
the standard presentation of $G\cG$ which
we suggestively called $K\cG^{(2)}$ before.  Thus
$\pi_1(K\cG) \simeq G\cG$.  It remains then
to show that the universal cover $\tKG$
is contractible.
\par
Since the universal cover is a 1-connected cubical complex, by
Gromov's theorem it suffices to show that the links of the vertices in
$\tKG$ are flag complexes.  To this end, let $v$ be a vertex in $\tKG$
and let $\cC$ be a complete graph contained in the 1-skeleton of the
link of $v$.  Each vertex of $\cC$ corresponds to an edge in
$\tKG$~---~the edge connecting the vertex to $v$~---~so the vertices of
$\cC$ correspond to the vertices of a subgraph $\cG_\cC$ of the
defining graph $\cG$.  In order to avoid confusion about which kind of
\lq\lq vertex" is being discussed, we use the phrase {\it link vertex}
for any vertex in $\cC$ which is contained in the link of $v$ and {\it
graph vertex} for any vertex in $\cG_\cC$.
\par
Because $\cC$ is a complete graph, each pair of link vertices has an
edge joining them.  Since this edge corresponds to a $2$-cell in
$\tKG$, the corresponding graph vertices must be joined by an edge.
Thus $\cG_\cC$ is also a complete graph (which is actually isomorphic
to $\cC$).  But by the construction of $\tKG$, $\cC$ must be embedded
in a cube in $\tKG$; hence $\cC$ is the $1$-skeleton of a simplex in
the link of $v$.
\hfill $\square$
\md
\vskip 1cm
\cl {\it {\rm 8.} Stability of higher invariants}
\bigskip
For many classes of groups it is known that
$\Sigma^1(G) = \Sigma^n(G)$ for all $n$.
For example, this stability of the BNS-invariants
is true for
$3$-manifold groups and $1$-relator groups \bs.
It is known that this
stability does not hold for arbitrary groups,
or even for arbitrary graph groups.
  For instance the direct product of two free
groups $F(a,b)$ and $F(c,d)$ is a graph group based
on a circuit of four edges.  If $\chi$ is the character sending
$a,b$ and $c$ to $1$ and $d$ to $0$,  then the kernel
of $\chi$ is finitely generated but not finitely presented
\stallings.  Thus by 1.1, $S(G\cG,Ker(\chi))\subseteq\Sigma^1(G\cG)$
while $S(G\cG,Ker(\chi))\not\subseteq\Sigma^2(G\cG)$,
so  $\Sigma^1(G\cG)\ne  \Sigma^2(G\cG)$.
There are, however, classes of graph groups
for which a certain amount of stability holds.
\par
Recall that to discuss the higher invariants
for graph groups
we must work with the maximal
subcomplex $\kplus$ of $\tKG$ contained in
${\tilde {\chi}}^{-1}([0,\infty))$.
\md
{\smps Definition.}  A simplicial graph is
{\it chordal} if and only if no full subgraph
generated by more than three vertices is a circuit.
\md
In \droms, C. Droms shows that graph groups based
on chordal graphs are coherent, i.e., every finitely generated subgroup
is finitely presented.  For a rational character
this  immediately
implies the following proposition, and for non-rational
characters the proposition can be proven using essentially
the same argument as appears in \droms.
\md
\proclaim Proposition 8.1. {\rm (Droms)} If $\cG$ is chordal, then
$\Sigma^1(G\cG) = \Sigma^2(G\cG)$.
\hfill $\square$
\par
Preliminary computations indicate that for a chordal graph $\cG$
it may be the case that
$\Sigma^1(G\cG) = \Sigma^k(G\cG)$ for all $k$, but we have
been unable to prove this.
\par
There are other situations where one has
this
amount of stability in the BNS-invariants.
The following lemma is proven using the ``$\Sigma^k$-Criteria"
in appendix B of \bs.  Since \bs\ is not widely available,
we write out a complete proof.
\par
\proclaim Lemma 8.2.  Let $\cG$ be a
simplicial graph with a vertex $v$
connected to every other vertex in $\cG$.
Let $x$ be the generator corresponding to $v$
and let $\chi$ be a character with $\chi(x) \not = 0$
(so $\cL(\chi)$ is connected and dominating).
Then $\chi$ defines a point in $\Sigma^k(G\cG)$
for all $k$.
\par
{\it Proof.} Let $v$ and $x$ be as in the statement
of the lemma.  Without loss of generality we may assume
$\chi(x)>0$.  For each integer $n$, let
$\tKG~{\buildrel f_n \over\longrightarrow}
~\tKG$ be defined by $f_n(p)=x^n\cdot p$, where $\cdot$
denotes the action of $G\cG$ on $\tKG$.  Extend this map
linearly to define for each real number $t$ a
 function $f_t$.  Notice for all $p\in\tKG$ and
for all $t$, ${\tilde{\chi}}(f_t(p))={\tilde{\chi}}(p) +
t\cdot {\chi}(x)$.  An
immediate consequence is that for all $p\in\tKG$, there exists $s\ge
0$ such that $f_t(p)\in\kplus$ for all $t\ge s$.
\par
We need to establish that $\kplus$ is $n$-connected for each $n$.
Since $\cL(\chi)$ is connected and dominating, $\kplus$ is connected
by 4.1.  Let $n\ge 1$ and $S^n~ {\buildrel \phi \over
\longrightarrow}~\kplus$ be an embedding of the $n$-sphere into
$\kplus$.  Since $\tKG$ is contractible, there is a map
$D^n~{\buildrel \Phi \over
\longrightarrow}~\tKG$ extending $\phi$ to a map of a disk into
$\tKG$.  Let $N(\Phi)$ be the minimum value of
$({\tilde{\chi}}\circ\Phi)(D^n)$.  While it is not necessarily true
that $\kplus ={\tilde \chi}^{-1} ([0,\infty))$, we may choose
$\epsilon >0$ such that ${\tilde
\chi}^{-1}([\epsilon, \infty))~\subseteq~\kplus$.  Since $\chi (x)>0$,
we may choose a constant
$K$ so that $K\chi(x) > N(\Phi) + \epsilon$.
\par
For $0\le t\le 1$, the family of maps
$\phi_t:S^n~\longrightarrow~\kplus$ given by
$\phi_t=f_{Kt}\circ\phi$ is a homotopy from $\phi$ to $f_K\circ
\phi$.
The map $f_K\circ\phi$ extends to the map
$f_K\circ \Phi~:~D^n~\longrightarrow~\tKG$;
 because the constant $K$ was chosen sufficiently large,
the image of $f_K\circ \Phi$ is actually contained in
$\kplus$.  It follows that $\kplus$ is
$n$-connected.
\hfill $\square$
\proclaim Theorem 8.3.  Let $G\cG$ be a graph group
which can be expressed
as the direct product of a non-trivial free abelian group and
a graph group with disconnected graph. Then
$\Sigma^1(G\cG) = \Sigma^k(G\cG)$ for all integers $k$.
\par
 {\it Proof.}  The underlying graph of such a graph group
decomposes into a complete graph $\cC$ corresponding to
the free abelian group, and pairwise disjoint connected graphs $\cG_i$,
with the only additional edges connecting each vertex of $\cC$ to
each vertex of $\cG_i$ for all $i$.
\par
Let $\chi$ represent a point in $\Sigma^1(G\cG)$.  By Theorem 4.1,
$\cL(\chi)$ is connected and dominating, so $\cL(\chi)$ must contain a
vertex of $\cC$.  By Lemma 8.2, $\chi$ represents a point in
$\Sigma^k(G\cG)$ for all $k$.  Since this is true for all characters
with finitely generated kernels and since $\Sigma^1(G) \supseteq
\Sigma^k(G)$ for any finitely generated group $G$, it follows that
$\Sigma^1(G\cG) = \Sigma^k(G\cG)$ for all $k$.
\hfill $\square$
\md
\vskip 1cm
\cl {\it References}
\bigskip
{\frenchspacing
\baselineskip = 12pt
\item {\bf 1.} {\smps R. Bieri, W.D. Neumann and R. Strebel}, `A
 Geometric Invariant of Discrete Groups', {\it Invent. Math.} 90
(1987) 451-477.\md
\item {\bf 2.} {\smps R. Bieri and B. Renz}, `Valuations on Free
Resolutions and Higher Geometric Invariants of Groups', {\it Comm.
Math. Helv.} 63 (1988) 464-497.\md
\item {\bf 3.} {\smps R. Bieri and R. Strebel}, `Geometric Invariants
for Discrete Groups', in preparation.\md
\item {\bf 4.} {\smps K.S. Brown}, {\it Cohomology of Groups}
(Springer-Verlag, New York, 1982).\md
\item {\bf 5.} {\smps P. Cartier and D. Foata}, `Probl\` emes
Combinatoires de Commutation et
R\'earrangements', Lecture Notes in
Math. 85 (Springer, Berlin, 1969).\md
\item {\bf 6.} {\smps C. Droms}, `Graph Groups, Coherence, and
Three-Manifolds', {\it J.  Algebra} 106  (1987) 484-489.\md
\item {\bf 7.} {\smps C. Droms}, `Subgroups of Graph Groups', {\it J.
Algebra} 110 (1987) 519--522.\md
\item {\bf 8.} {\smps D.B.A. Epstein, J.W. Cannon, D.F. Holt, M.S. Patterson,
W.P. Thurston}, {\it Word Processing in Groups} (Jones and Bartlett,
1992).\md
\item {\bf 9.} {\smps M. Gromov}, `Hyperbolic Groups', in {\it Essays
in Group Theory} ed. by S.M. Gersten, MSRI Publ. 8  (Springer-Verlag,
New York, 1987) 75-263.\md
\item {\bf 10.} {\smps S. Hermiller and J. Meier}, `Algorithms and
Geometry for Graph Products', to appear, {\it J.  Algebra}.\md
\item {\bf 11.} {\smps J. Meier}, `The topology of graph products
of groups', to appear, {\it Proc. Edin. Math. Soc.}.\md
\item {\bf 12.} {\smps H. Servatius, C. Droms, B. Servatius}, `Surface
Subgroups of Graph Groups', {\it Proc. Am. Math. Soc.} 106
(1989) 573--578.\md
\item {\bf 13.} {\smps J.R. Stallings}, `A Finitely Presented Group
Whose 3-Dimensional Integral Homology is not Finitely Generated',
{\it American J. of Math.} 85 (1963) 541-543.\md
\item {\bf 14.} {\smps L. VanWyk}, `Graph Groups are Biautomatic',
to appear, {\it J. Pure Appl. Alg.}.\md
\item {\bf 15.} {\smps C. Wrathall}, `The Word Problem for Free Partially
Commutative Groups', {\it J. Symbolic Computation}  6 (1988)
99-104.\md
}
\vskip 1cm
\settabs 8 \columns
\+& Department of Mathematics &&&& Department of Mathematics\cr
\+& Lafayette College                &&&& Binghamton University\cr
\+ & Easton                               &&&& Binghamton\cr
\+ & Pennsylvania   18042           &&&& New York 13902\cr
\+& USA                                   &&&& USA\cr
\md
\+& meierj@lafvax.lafayette.edu  &&&& vanwyk@math.binghamton.edu\cr
\vfill
\eject
\bye